\documentclass[a4 paper, 12pt]{amsart}
\usepackage[utf8]{inputenc}
\usepackage{amsmath}
\usepackage{amssymb}
\usepackage{amsfonts}

\usepackage{amsthm}

\numberwithin{equation}{section}

\newtheorem{theorem}{Theorem}
\newtheorem*{theorem*}{Theorem}

\newtheorem{corollary}{Corollary}
\newtheorem*{corollary*}{Corollary}

\newcommand{\eps}{\varepsilon}
\DeclareMathOperator{\spn}{span}

\usepackage{eufrak}
\newcommand{\ee}{\mathfrak{e}}

\DeclareRobustCommand{\divby}{%
	\mathrel{\text{\vbox{\baselineskip.65ex\lineskiplimit0pt\hbox{.}\hbox{.}\hbox{.}}}}%
}

\begin{document}
	\title[Bases in Hilbert and Banach spaces]{Almost Auerbach, Markushevich and Schauder bases in Hilbert and Banach spaces}
	
	\author{Anton Tselishchev}
	\address{Department of Mathematics, Bar-Ilan University, Ramat-Gan 5290002, Israel}
	\address{St. Petersburg Department, Steklov Math. Institute, Fontanka 27, St. Petersburg 191023 Russia}
	\email{celis-anton@yandex.ru}
	\thanks{Research supported by ISF Grant No.\ 1044/21 and by Russian Science Foundation grant no.\ 23-11-00171}
	
	\begin{abstract}
		For any sequence of positive numbers $(\eps_n)_{n=1}^\infty$ such that $\sum_{n=1}^\infty \eps_n = \infty$ we provide an explicit simple construction of $(1+\eps_n)$-bounded Markushevich basis in a separable Hilbert space which is not strong, or, in other terminology, is not hereditary complete; this condition on the sequence $(\eps_n)_{n=1}^\infty$ is sharp. Using a finite-dimensional version of this construction, Dvoretzky's theorem and a construction of Vershynin, we conclude that in any Banach space for any sequence of positive numbers $(\eps_n)_{n=1}^\infty$ such that $\sum_{n=1}^\infty \eps_n^2 = \infty$ there exists a $(1+\eps_n)$-bounded Markushevich basis which is not a Schauder basis after any permutation of its elements.
	\end{abstract}

\maketitle

\section{Introduction}

\subsection{The main definitions} The word ``basis'' has many different meanings in the Banach space theory. The questions of existence of bases in Banach spaces which possess some nice properties lie in the foundation of the geometry of Banach spaces, as well as other fields of Analysis. We begin with presenting certain definitions that will be used throughout the paper.

Let $X$ be a separable Banach space. A system of vectors $\{x_n\}_{n=1}^\infty$ in $X$ is \emph{complete} in $X$ if $\overline{\spn}\{x_n\}_{n=1}^\infty = X$. It is called \emph{minimal} if  $x_m\not\in \overline{\spn}\{x_n\}_{n\neq m}$ for every $m\ge 1$.

If $\{x_n\}_{n=1}^\infty$ is a complete and minimal system, then there exists a unique sequence of \emph{biorthogonal functionals} $\{x_n^*\}_{n=1}^\infty$ in $X^*$, i.e., $x_n^* (x_m) = \delta_{nm}$. In this case with any element $x\in X$ one can associate a formal series
\begin{equation}\label{Fourier_series}
	x\sim\sum_{n=1}^\infty x_n^*(x) x_n.
\end{equation}
If for every $x\in X$ this series converges to $x$, then $\{x_n\}_{n=1}^\infty$ is called a \emph{Schauder basis}. In this case there exists a constant $C$ such that for any $N>0$
\begin{equation}\label{projections_bounded}
	\Big\| \sum_{n=1}^N x_n^*(x)x_n \Big\| \le C\|x\|.
\end{equation}
The minimal constant such that this inequality holds is called the \emph{basis constant} of $(x_n)_{n=1}^\infty$.

It was a long-standing open problem whether every separable Banach space has a Schauder basis. It was solved in the negative by P. Enflo: in \cite{Enf73} he constructed the first example of a space without a Schauder basis.

It is natural therefore to consider some weaker notions of bases. It leads us to the next definition: a complete and minimal system $\{x_n\}_{n=1}^\infty$ in $X$ is called a \emph{Markushevich basis}, or \emph{M-basis}, if the sequence of its biorthogonal functionals $\{x_n^*\}$ is total, i.e., if $x_n^*(x) = 0$ for every $n\ge 1$, then $x=0$. The word ``basis'' here might be somewhat misleading: it does not mean that the series \eqref{Fourier_series} converges to $x$. Already in 1943 Markushevich in the paper \cite{Mar43} proved that every separable Banach space contains an M-basis.

An M-basis $\{x_n\}_{n=1}^\infty$ is called \emph{bounded} if there exists a constant $C$ such that $\|x_n\|\cdot \|x_n^*\|\le C$ for every $n\ge 1$. Since the inequality \eqref{projections_bounded} for a Schauder basis implies $\|x_n\|\cdot \|x_n^*\|\le 2C$, every Schauder basis is a bounded M-basis. Ovespian and Pełczy\'{n}ski in \cite{OP75} proved that every separable Banach space admits a bounded M-basis. Later in \cite{Pel76} Pełczy\'{n}ski proved that for every $\eps > 0$ any separable Banach space admits an M-basis $\{x_n\}_{n=1}^\infty$ with biorthogonal functionals $\{x_n^*\}_{n=1}^\infty$ such that $\|x_n\|\cdot \|x_n^*\|\le 1+\eps$.

An M-basis $\{x_n\}_{n=1}^\infty$ such that $\|x_n\| = \|x_n^*\| = 1$ is called an \emph{Auerbach basis}. Every finite dimensional Banach space contains an Auerbach basis --- this statement is known as Auerbach's lemma. However, the question whether every separable Banach space has an Auerbach basis is a long-standing open problem; in particular, it is not known for the space of continuous functions on a segment $C([0,1])$ (see e.g. \cite[Problem 113]{GMZ16}).

For a sequence of non-negative numbers $(\eps_n)_{n\ge 1}$ an M-basis $\{x_n\}_{n=1}^\infty$ is called \emph{$(1+\eps_n)$-bounded} if $\|x_n\|\cdot \|x_n^*\|\le 1+\eps_n$; if $\eps_n$ tends to $0$, then it is natural to call such M-bases \emph{almost Auerbach}. Vershynin in \cite{Ver99} proved that if $\sum_{n\ge 1} \eps_n^2 = \infty$, then every separable Banach space contains a $(1+\eps_n)$-bounded M-basis.

Finally, let us introduce one more definition. An M-basis $\{x_n\}_{n=1}^\infty$ is called \emph{strong} if $x\in\overline{\spn}\{x_n^*(x) x_n\}_{n=1}^\infty$ for any $x\in X$. If $X$ is a Hilbert space, then this property is equivalent to \emph{hereditary completeness} and it is related to so-called spectral synthesis of operators on Hilbert spaces, see e.g. \cite{Mar70}. The examples of bounded non-strong M-bases in a separable Hilbert space can be found e.g. in  \cite{Mar70, DN76, DNS77}. Obviously, if an M-basis is not strong, then it is not a Schauder basis after any permutation of its elements.

\subsection{Formulations of the results} An Auerbach basis in a Hilbert space is an orthonormal basis. It is therefore natural to ask if $(1+\eps_n)$-bounded M-bases also have some good properties for sufficiently small numbers $\eps_n$. In a recent paper \cite{RWZ24} it is shown that it is indeed the case: if $\sum_{n\ge 1} \eps_n < \infty$, then any (normalized) $(1+\eps_n)$-bounded M-basis in a separable Hilbert space must be a Riesz, or unconditional, basis (i.e., an image of an orthonormal basis under the action of a bounded invertible operator). Besides that, it was shown that the condition $\sum_{n\ge 1} \eps_n < \infty$ is sharp: if this series is divergent, then there exists a $(1+\eps_n)$-bounded M-basis which is not a Riesz basis.

Note that Riesz bases are always Schauder bases while the converse is not true (the examples of Schauder bases which are not Riesz bases in a Hilbert space can be found e.g. in \cite[Section 9.5]{AK16}). When it comes to the notion of a Schauder basis, in the paper \cite{RWZ24} it was shown that if the sequence $(\eps_n)_{n\ge 1}$ decays not very fast, then there exists a $(1+\eps_n)$-bounded M-basis which is not a Schauder basis after any permutation. However, the required assumption on this sequence is considerably stronger than the divergence of the series $\sum_{n\ge 1} \eps_n$. To be more precise, it was proved that if $\lim_{n\to\infty} n\eps_n = \infty$, then there exists a $(1+\eps_n)$-bounded M-basis in $\ell^2$ which is not a Schauder basis in any order. This basis was constructed from finite-dimensional ``blocks'' and hence it is a strong M-basis. Finally, using this finite-dimensional construction, the method of Vershynin from the paper \cite{Ver99} and Dvoretzky's theorem, the following corollary was proved in \cite{RWZ24}: if $\lim_{n\to\infty} n\eps_n = \infty$ and $\sum_{n=1}^\infty\eps_n^2 = \infty$, then in any separable Banach space $X$ there exists a $(1+\eps_n)$-bounded M-basis which cannot be ordered to become a Schauder basis.

Our goal is to improve the aforementioned results from \cite{RWZ24}. In particular, we will get rid of the condition $\lim_{n\to \infty} n\eps_n = \infty$. At first, we will present a simple explicit construction which proves the following statement.

\begin{theorem}
	\label{example}
	For any sequence of nonnegative numbers $(\eps_n)_{n\ge 1}$ such that $\sum_{n=1}^\infty \eps_n = \infty$ there exists a $(1+\eps_n)$-bounded M-basis in $\ell^2$ which is not strong.
\end{theorem}

Non-strong M-bases are not Riesz bases and moreover, as we mentioned before, they do not form Schauder bases after any permutation.

In order to obtain a result for arbitrary Banach spaces using Dvoretzky's theorem, we will need the following finite-dimensional version of our construction.

\begin{theorem}
	Suppose that $(\eps_n)_{n\ge 1}$ is a sequence of nonnegative numbers such that $\sum_{n=1}^\infty \eps_n = \infty$. Then for any $C > 0$ there exists an integer number $N$, a finite-dimensional Hilbert space $H$, $\dim H = N$, and a biorthogonal system $\{(x_n, x_n^*)\}_{n=1}^N$ in $H$ such that $\|x_n\|\cdot \|x_n^*\|\le 1+\eps_n$ and for any bijection $\pi:\{1,\ldots, N\}\to \{1,\ldots, N\}$ the basis constant of $(x_{\pi(k)})_{k=1}^N$ is at least $C$.
\end{theorem}

Now, using Vershynin's construction from the paper \cite{Ver99} and Dvoretzky's theorem, the following corollary can be proved in exactly the same way as it is done in \cite{RWZ24} --- but the extra condition $\lim_{n\to\infty} n\eps_n = \infty$ is no longer required (note that if $\sum_{n=1}^\infty\eps_n^2 = \infty$, then of course also $\sum_{n=1}^\infty\eps_n = \infty$).

\begin{corollary}
	Suppose that $(\eps_n)_{n\ge 1}$ is a sequence of nonnegative numbers such that $\sum_{n=1}^\infty\eps_n^2 = \infty$. Then in any separable Banach space $X$ there exists a $(1+\eps_n)$-bounded M-basis which cannot be ordered to become a Schauder basis.	
\end{corollary}

Besides that, we can use Theorem~2 in order to construct a $(1+\eps_n)$-bounded strong M-basis which is not a Schauder basis after any permutation of its elements:

\begin{corollary}
	For any sequence of nonnegative numbers $(\eps_n)_{n\ge 1}$ such that $\sum_{n=1}^\infty \eps_n = \infty$ there exists a $(1+\eps_n)$-bounded strong M-basis in $\ell^2$ which is not a Schauder basis after any permutation. 
\end{corollary}

Indeed, the construction is straightforward. We can take an arbitrary sequence $C_m$ tending to infinity. Then, using Theorem~2, at first we find the space $H_1$ with $\dim H_1 = N_1$ and a $(1+\eps_n)_{n=1}^{N_1}$-bounded biorthogonal system in $H_1$ such that after any permutation its basis constant is at least $C_1$. After that we use Theorem~2 once again for a sequence $(\eps_{n+N_1})_{n\ge 1}$ and find a space $H_2$ with $\dim H_2 = N_2$ and a $(1+\eps_{n+N_1})_{n=1}^{N_2}$-bounded biorthogonal system in $H_2$ such that after any permutation its basis constant is at least $C_2$. We continue this process, and then we take the union of the constructed systems in $\bigoplus_{m=1}^\infty H_m = \ell^2$. Since our system in constructed from finite-dimensional ``blocks'', it is easy to see that it forms a strong M-basis.

\section{Proof of Theorem 1}

Let us prove Theorem~\ref{example}. First of all, we can assume that all numbers $\eps_n$ are strictly positive (because we can exclude all zeroes from the sequence $(\eps_n)_{n\ge 1}$, construct the required M-basis and then add extra coordinates to our Hilbert space and extra orthonormal vectors to our construction). After that, replacing, if necessary, each number $\eps_n$ with a smaller number we can assume that $\eps_n \to 0$ (but still the series $\sum_{n\ge 1} \eps_n$ is divergent). Finally, after a sufficient permutation, we can also assume that the sequence $(\eps_n)$ is non-increasing: $\eps_{n+1}\le \eps_n$ for every $n$.

Now let us define by induction the following sequence of positive numbers $(c_n)_{n=1}^\infty$:

\begin{equation}
	c_1=\eps_1; \quad c_{n+1}=\frac{\eps_{n+1}}{1+c_1+c_2+ \ldots + c_n}.
\end{equation}

We claim that $\sum_{n=1}^\infty c_n = \infty$. Indeed, suppose that $\sum_{n=1}^\infty c_n = A < \infty$. Then
\begin{equation*}
	c_{n+1} = \frac{\eps_{n+1}}{1+c_1+\ldots + c_n}\ge \frac{\eps_{n+1}}{1 + A}
\end{equation*}
and therefore
\begin{equation*}
	\sum_{k=1}^\infty \eps_k \le (1+A)\sum_{k=1}^\infty c_k = A(1+A) < \infty,
\end{equation*}
which contradicts the assumption of the divergence of series $\sum_{n\ge 1}\eps_n$.

It is not difficult to see that the sequence $(c_n)_{n=1}^\infty$ is non-increasing: for every $n\ge 1$ we have
\begin{equation}\label{non-increase}
	c_{n+1} = \frac{\eps_{n+1}}{1+c_1+\ldots + c_n} \le \frac{\eps_n}{1+c_1+\ldots + c_n}\le \frac{\eps_n}{1+c_1+\ldots + c_{n-1}} = c_n.
\end{equation}

Since the sequence $(c_n)$ is non-increasing and the series $\sum_{n=1}^\infty c_n$ diverges, we also have
\begin{align}
	\sum_{k=1}^\infty c_{2k} &= \infty;\label{div-even}\\
	\sum_{l=1}^\infty c_{2l-1} &= \infty.\label{div-odd}
\end{align}

Let us define the following numbers for $m,n \ge 1$:
\begin{align}
	\beta_m &= c_m^{1/2}; \label{betadef}\\
	\alpha_{m,n} &= \begin{cases}
		-\beta_m\beta_n, &m>n\\
		0, &m\le n
	\end{cases}. \label{alphadef}
\end{align}

Suppose now that $\{\ee\} \cup \{e_n\}_{n=1}^\infty$ is an orthonormal basis in $\ell^2$. We define the vectors $\{x_n\}_{n=1}^\infty$ and $\{y_n\}_{n=1}^\infty$ in our space $\ell^2$ as follows:
\begin{align}
	x_{2l-1} &= e_{2l-1}; \  y_{2l-1} = \beta_{2l-1}\ee + e_{2l-1} + \sum_{k=1}^\infty \alpha_{2l-1,2k}e_{2k},\ l\ge 1 \label{def1}\\
	y_{2k}&=e_{2k}; \qquad x_{2k}=\beta_{2k}\ee + e_{2k} + \sum_{l=1}^\infty \alpha_{2k, 2l-1}e_{2l-1},\ \  k\ge 1.\label{def2}
\end{align}
Note that each sum in this definition contains only a finite number of nonzero summands. We claim that the system $\{y_n\}_{n=1}^\infty$ is biorthogonal to $\{x_n\}_{n=1}^\infty$, both these systems are complete in $\ell^2$ but $\{x_n\}_{n=1}^\infty$ is not a strong M-basis and $\|x_n\|\cdot \|y_n\|\le 1 + \eps_n$ for every $n\ge 1$. These properties imply that $\{x_n\}_{n=1}^\infty$ serves as an example for Theorem~\ref{example}. In what follows we will verify these properties.

\subsection{Completeness}
	
	Let us prove that the system $\{x_n\}_{n=1}^\infty$ is complete in $\ell^2$ (the completeness of $\{y_n\}_{n=1}^\infty$ can be proved in the same way). Suppose that $h\in\ell^2$ and $\langle h, x_n \rangle = 0$ for every $n\ge 1$. Then in particular $\langle h, x_{2l-1}\rangle = \langle h, e_{2l-1} \rangle = 0$. Therefore, the relation $\langle h, x_{2k} \rangle = 0$ is equivalent to $\langle h, e_{2k}\rangle = -\beta_{2k} \langle h, \ee \rangle$. If we now denote $c = \langle h, \ee \rangle$, we see that
	\begin{equation*}
		\sum_{k=1}^\infty \langle h, e_{2k} \rangle^2 = c^2 \sum_{k=1}^\infty \beta_{2k}^2 = c^2 \sum_{k=1}^\infty c_{2k}.
	\end{equation*}
	The latter series diverges (see formula \eqref{div-even}), hence, since $h\in\ell^2$ we infer that $c$ must be equal to $0$ --- but then $\langle h, \ee \rangle = \langle h, e_n \rangle = 0$ for every $n\ge 1$ which means that $h=0$ in $\ell^2$. This proves the completeness of the system $\{x_n\}_{n=1}^\infty$ and, as we mentioned, completeness of the system $\{y_n\}_{n=1}^\infty$ follows in a similar way.
	
	\subsection{Biorthogonality}
	The relations $\langle x_{2l-1}, y_{2k}\rangle = 0$,  $\langle x_{2l-1}, y_{2m-1} \rangle = \delta_{lm}$ and $\langle x_{2k}, y_{2m} \rangle = \delta_{km}$ are trivial. We only need to check that $\langle x_{2k}, y_{2l-1}\rangle = 0$. But the direct computation shows that $\langle x_{2k}, y_{2l-1} \rangle = \beta_{2l-1}\beta_{2k} + \alpha_{2l-1, 2k} + \alpha_{2k, 2l-1}$. This quantity equals to $0$ by the definition of the numbers $\alpha_{m,n}$ (see formula \eqref{alphadef}).
	
	\subsection{The basis is not strong}
	We have just shown that $\{x_n\}_{n=1}^\infty$ is an M-basis in $\ell^2$ and $\{y_n\}_{n=1}^\infty$ is its biorthogonal system. However, it is straightforward to see that this M-basis is not strong. Indeed $\langle \ee, y_n \rangle \neq 0$ only for odd numbers $n$ and hence
	\begin{equation*}
		\ee \not\in \overline{\spn}\{e_{2l-1}\}_{l=1}^\infty = \overline{\spn}\{\langle \ee,  y_n \rangle x_n\}_{n=1}^\infty.
	\end{equation*}

	\subsection{The norm estimate}
	Finally, it remains only to check that $\|x_n\| \cdot \|y_n\| \le 1+\eps_n$ for all $n\ge 1$. Let us check it for even numbers $n=2k$ (the computation for odd numbers is completely similar). In this case $\|y_{2k}\| = 1$ and
	\begin{multline*}
		\|x_{2k}\|^2 = 1 + \beta_{2k}^2 + \sum_{l=1}^\infty \alpha_{2k, 2l-1}^2 = 1+\beta_{2k}^2+\sum_{l=1}^k \beta_{2k}^2 \beta_{2l-1}^2 = 1+c_{2k}\Big( 1 +  \sum_{l=1}^k c_{2l-1}\Big)\\
		\le 1+c_{2k}\Big(1+\sum_{n=1}^{2k-1} c_n\Big)=1+\eps_{2k}.
	\end{multline*}
	Using this inequality, we get that $\|x_{2k}\| \le (1+\eps_{2k})^{1/2} \le 1+\eps_{2k}$, and we are done.
	
	\section{Proof of Theorem 2}
	
	As in the proof of Theorem~1, we can assume that the numbers $\eps_n$ are positive and the sequence $(\eps_n)_{n\ge 1}$ tends to $0$. Besides, we can assume that  $\eps_{n+1}\le \eps_n$ for all $n\ge 1$. Indeed, note that if we construct a Hilbert space $H$ with a basis $(x_n)_{n=1}^N$ in it which satisfies the conditions of Theorem~2 for non-increasing permutation of the numbers $(\eps_n)_{n\ge 1}$, then, for a given arbitrary order of $(\eps_n)_{n\ge 1}$, we can consider the Hilbert space $\widetilde{H}$ of sufficiently large dimension $K$, fix an orthonormal basis $\{w_n\}_{n=1}^K$  in it, and take the space $\widetilde{H}\oplus H$ with a sufficient permutation of the basis $\{w_n\oplus 0\}_{n=1}^K\cup \{0\oplus x_n\}_{n=1}^N$ in it.  
	
	Let us define the numbers $\beta_n$ and $\alpha_{m,n}$ in the same way as in the proof of Theorem~1 (see formulas \eqref{betadef} and \eqref{alphadef}). Recall that, due to our assumptions, as before, the sequence $(\beta_n)$ is non-increasing (see the formulas \eqref{non-increase} and \eqref{betadef}). The formulas \eqref{div-even} and \eqref{div-odd} can be rewritten as follows:
	\begin{align}
		\sum_{k=1}^\infty \beta_{2k}^2 &= \infty;\label{div-even_beta}\\
		\sum_{l=1}^\infty \beta_{2l-1}^2 &= \infty.\label{div-odd_beta}
	\end{align}
	Hence we can fix a number $L$ such that
	\begin{equation}\label{L_is_big}
		\sum_{l=1}^L \beta_{2l-1}^2 \ge 64C^2.
	\end{equation}

	After that we can find a number $M > L$ such that
	\begin{equation}\label{M_is_large}
		\sum_{k=1}^L \beta_{2k}^2 \le \frac{1}{8}\sum_{k=1}^M \beta_{2k}^2.
	\end{equation}
	We can also assume that
	\begin{equation}\label{A_is_big}
	\sum_{k=1}^M \beta_{2k}^2 > 1
	\end{equation}
	We put $N=2M$.
	
	Consider, as before, the space $\ell^2$ with an orthonormal basis $\{\ee\} \cup \{e_n\}_{n=1}^\infty$. Recall the formulas \eqref{def1} and \eqref{def2} which define the biorthogonal system $\{(x_n, y_n)\}_{n=1}^\infty$ in $\ell^2$. For the reader's convenience, we rewrite them here is the following form:
	\begin{align}
		x_{2l-1} &= e_{2l-1}; \ y_{2l-1}=\beta_{2l-1}\ee + e_{2l-1} - \sum_{k=1}^{l-1} (\beta_{2l-1}\beta_{2k}) e_{2k}, \  l\ge 1;\label{defodd}\\
		y_{2k}&=e_{2k}; \ \ \ \ \ x_{2k} = \beta_{2k}\ee + e_{2k}-\sum_{l=1}^{k}(\beta_{2k}\beta_{2l-1})e_{2l-1}, \ \  k\ge 1.\label{defeven}
	\end{align}
	
	 Now we take $H=\mathrm{span}\{x_1, x_2, \ldots, x_{2M}\}\subset\ell^2$ and denote by $P$ an orthogonal projection onto $H$ in $\ell^2$. Then, if we put $x_n^* = Py_n$, we get that the system $(x_n^*)_{n=1}^{2M}$ is biorthogonal to $(x_n)_{n\ge 1}$ in $H$. Besides that, in the proof of Theorem~1 we checked that $\|x_n\|\cdot \|y_n\|\le 1+\eps_n$ for every $n$, and, since orthogonal projection does not increase the norm of the vector, $\|x_n\|\cdot \|x_n^*\|\le 1+\eps_n$. It remains to check that for any permutation $\pi:\{1,\ldots, 2M\}\to \{1,\ldots, 2M\}$ the basis constant of $(x_{\pi(n)})_{n=1}^{2M}$ is at least $C$.
	
	We denote 
	$$
	B=\sum_{k=1}^M \beta_{2k}^2.
	$$
	Note that $B > 1$ by \eqref{A_is_big}. Now consider the following coefficients:
	\begin{align}
		\gamma_{2k}&=\frac{\beta_{2k}}{B}, \ 1\le k\le M;\label{defeven-gamma}\\
		\gamma_{2l-1}&=\frac{\beta_{2l-1}}{B}\Big( \sum_{k=l}^M \beta_{2k}^2 \Big),\  1\le l \le M.\label{defodd-gamma}
	\end{align}
	Finally, let us consider the following element of $H$:
	\begin{equation}
		z=\sum_{n=1}^{2M} \gamma_n x_n.
	\end{equation}

Let us compute $\|z\|$. Obviously, $z\in\mathrm{span}\{\ee, e_1, e_2,\ldots e_{2M}\}$ in $\ell^2$. A direct computation shows that for $1\le k, l \le M$
\begin{equation*}
	\langle z, \ee \rangle = 1; \ \langle z, e_{2l-1} \rangle = 0; \ \langle z, e_{2k}\rangle = \frac{\beta_{2k}}{B}.
\end{equation*}
Hence
$$
\|z\|=\Big( 1+\sum_{k=1}^M \frac{\beta_{2k}^2}{B^2} \Big)^{1/2} = (1+1/B)^{1/2}\le 2.
$$

Let us now fix a permutation $\pi:\{1,\ldots, 2M\}\to \{1,\ldots, 2M\}$. If we find a number $t$ such that
$$
\Big\| \sum_{m=1}^t \gamma_{\pi(m)}x_{\pi(m)} \Big\| \ge 2C,
$$
then we are done.

There exists a unique number $t$ such that
\begin{equation}\label{t-def}
	\sum_{\{m\le t :\,  \pi(m)\divby 2\}} \beta_{\pi(m)}^2 \ge \frac{B}{2}, \quad \sum_{\{m\le t-1 :\,  \pi(m)\divby 2\}} \beta_{\pi(m)}^2 < \frac{B}{2}.
\end{equation}
Since for every $k\ge 1$ we have $\beta_2\ge\beta_{2k}$, these relations imply
\begin{equation}\label{est_above}
	\sum_{\{m\le t :\,  \pi(m)\divby 2\}} \beta_{\pi(m)}^2 < \frac{B}{2} + \beta_{\pi(t)}^2\le \frac{B}{2} + \beta_2^2.
\end{equation}
We put
$$
w=\sum_{m=1}^t \gamma_{\pi(m)} x_{\pi(m)}\in H\subset\ell^2.
$$
We will prove now that $\|w\|\ge 2C$.

Let us compute the quantity $\langle w, e_{2l-1} \rangle$ for $l=1, 2,\ldots, L$. Note that $\langle x_n, e_{2l-1}\rangle\neq 0 $ if either $n$ is even and $n \ge 2l $ or $n=2l-1$. Hence, if $2l-1\not\in \{\pi(1), \pi(2), \ldots, \pi(t)\}$, then, using the formulas \eqref{defeven} and \eqref{defeven-gamma}, we get
\begin{equation}\label{firstcase}
	\langle w, e_{2l-1} \rangle = - \frac{1}{B} \sum_{\{m\le t :\,  \pi(m)\divby 2, \, \pi(m)\ge 2l\}} \beta_{\pi(m)}^2 \beta_{2l-1}.
\end{equation}
In this case, using the formulas \eqref{M_is_large} and \eqref{t-def}, we can estimate this quantity as follows:
\begin{multline}\label{est_first_case}
	|\langle w, e_{2l-1}\rangle| \ge \frac{\beta_{2l-1}}{B}\Big( \sum_{\{m\le t :\,  \pi(m)\divby 2\}} \beta_{\pi(m)}^2  - \sum_{k=1}^{l-1} \beta_{2k}^2 \Big) \\ \ge \frac{\beta_{2l-1}}{B}\Big( \sum_{\{m\le t :\,  \pi(m)\divby 2\}} \beta_{\pi(m)}^2  - \sum_{k=1}^{L} \beta_{2k}^2 \Big)  \ge \frac{\beta_{2l-1}}{B} \Big(\frac{B}{2} - \frac{B}{8}\Big) \ge \frac{\beta_{2l-1}}{4}.
\end{multline}

Now we will compute  $\langle w, e_{2l-1} \rangle$ in the case when $2l-1 \in \{\pi(1), \pi(2), \ldots, \pi(t)\}$. In this case we should add $\gamma_{2l-1}$ to the right-hand side of the formula \eqref{firstcase} and get
\begin{multline*}
		\langle w, e_{2l-1} \rangle = - \frac{1}{B} \sum_{\{m\le t :\,  \pi(m)\divby 2, \, \pi(m)\ge 2l\}} \beta_{\pi(m)}^2 \beta_{2l-1} + \frac{\beta_{2l-1}}{B}\Big( \sum_{k=l}^M \beta_{2k}^2 \Big) \\ = \frac{\beta_{2l-1}}{B} \sum_{\{m > t :\,  \pi(m)\divby 2, \, \pi(m)\ge 2l\}}\beta_{\pi(m)}^2.
\end{multline*}
Using the formulas \eqref{est_above} and \eqref{M_is_large}, in this case we can write:
\begin{multline}\label{est_second_case}
	|\langle w, e_{2l-1} \rangle| \ge \frac{\beta_{2l-1}}{B}\Big( B - \sum_{\{m\le t :\,  \pi(m)\divby 2\}} \beta_{\pi(m)}^2  - \sum_{k=1}^{{l-1}} \beta_{2k}^2  \Big) \\ \ge \frac{\beta_{2l-1}}{B}\Big( B - \frac{B}{2} - \beta_2^2 - \sum_{k=1}^{L} \beta_{2k}^2  \Big)\ge \frac{\beta_{2l-1}}{B} \Big( \frac{B}{2} - \frac{B}{8} - \frac{B}{8} \Big) = \frac{\beta_{2l-1}}{4}.
\end{multline}

Summing up, the formulas \eqref{est_first_case} and \eqref{est_second_case} imply that in either case $|\langle w, e_{2l-1}\rangle| \ge \frac{\beta_{2l-1}}{4}$ for $l=1, 2, \ldots L$. Hence, by \eqref{L_is_big},
\begin{equation*}
	\| w \|\ge \Big( \sum_{l=1}^L \langle w, e_{2l-1} \rangle^2 \Big)^{1/2} \ge \frac{1}{4} \Big( \sum_{l=1}^L \beta_{2l-1}^2 \Big)^{1/2} \ge 2C,
\end{equation*}
and the theorem is proved.

\section*{Acknowledgements} The author is grateful to Nir Lev for fruitful discussions.

\end{document}